\documentclass{amsart}
\usepackage{amssymb}

\begin{document}

\title{Euler and the strong law of small numbers}

\author{Karl Dilcher}
\address{Department of Mathematics and Statistics\\
         Dalhousie University\\
         Halifax, Nova Scotia, B3H 4R2, Canada}
\email{dilcher@mathstat.dal.ca}
\author{Christophe Vignat}
\address{LSS-Supelec, Universit\'e Paris-Sud, Orsay, France and Department of
Mathematics, Tulane University, New Orleans, LA 70118, USA}
\email{cvignat@tulane.edu}

\setcounter{equation}{0}

\maketitle

In addition to modern computer algebra systems and bibliographic databases,
old-fashioned tables of integrals, sums and products remain important research
tools in mathematics. This is especially the case for the present authors 
whose research is in the classical areas of analysis, combinatorics, and 
number theory.

What follows is a cautionary tale about the use of tables, along with a 
reminder, as much to ourselves as to the reader, to go back to the sources
if at all possible. But before we continue, an important caution to the
reader of this note: {\it Some of the identities that follow are incorrect, 
while others may be misleading\/}. Also, some notations will be inconsistent, 
to preserve the historical context.

In the process of looking up some infinite products in the excellent and
very useful handbook by E.~R.~Hansen \cite{Ha}, we came across
a curious identity for the Euler numbers $E_n$. These numbers can be defined
by the exponential generating function
\begin{equation}\label{1}
\frac{2}{e^t+e^{-t}} = \sum_{n=0}^\infty E_n\frac{t^n}{n!}.
\end{equation}
Since the left-hand side of \eqref{1} is an even function, we have 
$E_{2k+1}=0$ for all $k\geq 0$. Furthermore, it can be shown that the 
even-index Euler numbers are integers, and $E_0=1$, $E_2=-1$, $E_4=5$,
$E_6=-61$, $E_8=1385,\ldots$. For further properties, see, e.g., \cite{AS} or
its successor volume \cite{DLMF}. The identity in question is (89.8.3) on
p.~486 of \cite{Ha}, namely
\begin{equation}\label{2}
\prod_{k=1}^\infty\left[1-\frac{(-1)^k}{(2k+1)^{2n+1}}\right]
= \frac{1}{2(2n)!}\left|E_{2n}\right|\left(\frac{\pi}{2}\right)^{2n+1}.
\end{equation}

A quick numerical check with {\it Maple} and {\it Mathematica} for a few 
small values of $n$ revealed that the identity is definitely 
{\it incorrect\/}. Before we realized what was going on, we
found the identity, in the same form, in the first volume of the well-known
multi-volume tables by Prudnikov et al. \cite[p.~754]{PrE}. No reference is 
given, and the identity also appears in the Russian original \cite{PrR}, again
in the same form.

One of the useful features of Hansen's handbook \cite{Ha} is the fact that most
entries come with one or more references to the literature, or at least to 
other entries in the table. The reference associated to identity \eqref{2}, 
i.e., Hansen's (89.8.3), is given as identity (1131) in the smaller and older 
collection of formulas by Jolley \cite{Jo}. The identity in question is listed 
on pp.~238--239 as
\begin{equation}\label{3}
2(2n)!\left(\frac{2}{\pi}\right)^{2n+1}\left(1+\frac{1}{3^{2n+1}}\right)
\left(1-\frac{1}{5^{2n+1}}\right)\left(1+\frac{1}{7^{2n+1}}\right)\ldots
= E_n^* = E_{2n},
\end{equation}
where $E_n^*$ and $E_{2n}$ are notations of the Euler numbers that differ from
the almost universally accepted modern usage in \eqref{1} and \eqref{2}, namely
$E_0^*=E_1^*=1, E_2^*=5, E_3^*=61, E_4^*=1385,\ldots$.

Jolley also provides references, and for his identity (1131), i.e., \eqref{3}
above, the reader is referred to p.~365 of the famous old Algebra book by
Chrystal \cite{Ch}; it is actually part II of this two-volume work, a fact 
that is not mentioned in \cite{Jo}. The identity in question turns out to be
(15) on p.~365 of \cite{Ch}, namely
\begin{equation}\label{4}
E_m=2(2m)!\left(\frac{2}{\pi}\right)^{2m+1}\left/
\left(1+\frac{1}{3^{2m+1}}\right)
\left(1-\frac{1}{5^{2m+1}}\right)\left(1+\frac{1}{7^{2m+1}}\right)\ldots,\right.
\end{equation}
where yet another notation for the Euler numbers is used, namely
$E_1=1, E_2=5, E_3=61, E_4=1385,\ldots$ (see \cite[p.~342]{Ch}). Note that the
large fraction sign in \eqref{4} is missing in \eqref{3}.

In addition to hinting at a proof of this last identity, Chrystal refers the
reader to Euler by writing, ``See again Euler, {\it Introd. in Anal. Inf.\/},
\S 284" in a footnote on p.~365. Euler's influential book \cite{Eu1L} is 
available in English translation \cite{Eu1E}, and \S 284 can be found on
pp.~239--241 of \cite{Eu1L}, or on pp.~244--245 of \cite{Eu1E}, where we find
\begin{equation}\label{5}
A = 1-\frac{1}{3^n}+\frac{1}{5^n}-\frac{1}{7^n}+\frac{1}{9^n}-\frac{1}{11^n}
+\frac{1}{13^n}-\frac{1}{15^n}+\&c.
\end{equation}
and at the end of \S 284, Euler writes,
\begin{equation}\label{6}
A = \frac{3^n}{3^n+1}\cdot\frac{5^n}{5^n-1}\cdot\frac{7^n}{7^n+1}\cdot
\frac{11^n}{11^n+1}\cdot\frac{13^n}{13^n-1}\cdot\frac{17^n}{17^n-1}\cdot\&c,
\end{equation}
{\it where the powers of all the prime numbers in the numerators occur and
the denominators are increased or decreased by $1$ depending on whether the
number has the form $4m-1$ or $4m+1$.} \cite[p.~245]{Eu1E}.

The identity \eqref{6} already gives us a glimpse of what might have gone 
wrong on the way towards the infinite product \eqref{2}. But first we need
to establish a connection between Euler's expressions \eqref{5} and \eqref{6}
and what were later called the Euler numbers.

In another well-known book \cite{Eu2L}, published a few years after 
\cite{Eu1L}, Euler expressed the numbers $A$ in \eqref{5} in terms of the
Taylor coefficients of the secant function, which by \eqref{1} are closely
related to the Euler numbers. For instance, one of the explicit identities
that can be found in \S 224 in \cite[p.~542]{Eu2L}, or in German translation
in \cite[p.~259]{Eu2G}, is
\begin{equation}\label{7}
1-\frac{1}{3^9}+\frac{1}{5^9}-\frac{1}{7^9}+\frac{1}{9^9}-\&c.
=\frac{\varepsilon}{1\cdot 2\ldots 8}\cdot\frac{\pi^9}{2^{10}},
\end{equation}
where $\varepsilon=1385$. Since this note deals with incorrect identities,
we must mention in passing that there are typographical errors in the original
identity (and neighboring ones) on p.~542 of \cite{Eu2L}; however, only two
pages further they are correctly printed; it is also correct in \cite{Eu2G}.

The rest of the story is now easy to piece together. The numbers 
$\alpha, \beta, \ldots, \varepsilon, \ldots$ used by Euler in 
\cite[p.~542]{Eu2L} correspond to $1, E_1, \ldots, E_4, \ldots$ in Chrystal's
notation, and \eqref{4} is indeed the general form of Euler's identity ---
if the sequence $3, 5, 7, \ldots$ is interpreted as the beginning of the 
sequence of odd primes, rather than the sequence of odd integers greater than
1. Therefore, the origin of the incorrect formula \eqref{2} is quite likely
Chrystal's identity \eqref{4}. Showing just one more term, along the lines of
Euler's identity \eqref{6}, would have avoided all this. However, whether or
not Jolley misinterpreted the sequence $3, 5, 7, \ldots$, the identity (3)
still contains the mistake of the missing fraction sign which then made it
into Hansen's identity (2).

In the end, we shouldn't blame Chrystal too much. Given that his book is written
in great detail, even a moderately attentive reader would realize that the 
product in \eqref{4} had to be over the odd primes. This is all the more so
as Chrystal remarks that a previous identity is transformed into his (15),
i.e., identity \eqref{4} above, ``{\it in the same way as before}." He
apparently refers to identity (8) in \cite[p.~364]{Ch}, namely
\begin{equation}\label{8}
B_m=2(2m)!\left/(2\pi)^{2m}\left(1-1/2^{2m}\right)
\left(1-1/3^{2m}\right)\left(1-1/5^{2m}\right)\ldots\right.
\end{equation}
(once again in Chrystal's notation), where $B_m$ is the $m$th Bernoulli 
number in the historical notation that has $B_1=1/6, B_2=1/30, B_3=1/42,
B_4=1/30, B_5=5/66,\ldots$; for different notations see, e.g.,
\cite[Ch.~24]{DLMF}. This last identity \eqref{8} is closely related to the
Euler product for the Riemann zeta function, especially if we compare it with
the following famous formula named after Euler:
\begin{equation}\label{9}
B_m=\frac{2(2m)!}{(2\pi)^{2m}}
\left\{\frac{1}{1^{2m}}+\frac{1}{2^{2m}}+\frac{1}{3^{2m}}+\ldots\right\},
\end{equation} 
again reproduced as in \cite[p.~363]{Ch}. While there is much less danger of
the product in \eqref{8} to be misunderstood, Euler himself showed more 
terms in the analogs of \eqref{8} and \eqref{9}; see \S 283 in \cite{Eu1L} 
or \cite{Eu1E}.

We already mentioned that the product in \eqref{8} is, essentially, the Euler
product for the Riemann zeta function. Similarly, the product in \eqref{4} is
the Euler product of an appropriate $L$-series (in fact, the series \eqref{5}),
and both are special cases of Euler products of Dirichlet $L$-series; see,
e.g., \cite[p.~162ff.]{Co}.

Before we close, let us reiterate that some of the identities in this note
are incorrect or misleading. Indeed, the reader will have realized that
\eqref{2} and \eqref{3} are {\it false}, and \eqref{4} is correct only when
interpreted as a product over the odd primes. Along with the incorrect entry
(89.8.3) in \cite{Ha}, i.e., \eqref{2} above, two consequences are mentioned, 
namely (89.4.11) for $n=0$, and (89.6.12) for $n=1$. The first one of these
is correct since it is in a somewhat different form, taken from an identity 
in the well-known classical book by Bromwich \cite[p.~224, Ex.~9]{Br}.
The special case of (89.4.11), namely $w=-1$, that is relevant here  is
\begin{equation}\label{10}
\prod_{k=1}^\infty\left[1-\frac{(-1)^k}{2k+1}\right]
= \frac{\pi\sqrt{2}}{4}.
\end{equation}
However, the second identity, (89.6.12), is indeed false, the corrected 
version being
\begin{equation}\label{11}
\prod_{k=1}^\infty\left[1-\frac{(-1)^k}{(2k+1)^3}\right]
=\frac{\pi}{12} + \frac{\pi\sqrt{2}}{12}\cosh\left(\frac{\pi\sqrt{2}}{4}\right).
\end{equation}
In an attempt to prove \eqref{11}, we tried to use Hansen's formula (89.6.2) 
which, interestingly, also turned out to be  incorrect. In this case, however,
there was only a small
typographical error (it should read $t=e^{\pi i/3}$). We first obtained the 
correct version \eqref{11} symbolically with the computer algebra system
{\it Mathematica\/}. We were then able to prove its extension to arbitrary 
(even or odd) powers of the denominators $2k+1$, using properties of the 
Gamma function. But this is another story.

The mistake that has been the main subject of this note is a good example of 
the Strong Law of Small Numbers which Richard K.~Guy \cite{Gu} formulated as 
follows: ``{\it There aren't enough small numbers to meet the many demands
made of them.}" One of the corollaries stated in \cite{Gu} captures the 
situation even better: ``{\it Superficial similarities spawn spurious 
statements.}"

However, we would like to give the last word to Pierre-Simon Laplace and his
famous dictum: ``{\it Lisez Euler, lisez Euler, c'est notre ma{\^ i}tre \`a
tous.}" (``{\it Read Euler, read Euler, he is the master of us all.}") 
\cite{Wi}. And to do so, the best place to begin is the wonderful Euler 
Archive \cite{EA}.

\end{document}